\documentclass[12pt,twoside]{amsart}
 
\usepackage{amsmath,comment} 
\usepackage{amssymb} 

\numberwithin{equation}{section}
  
\newtheorem{thm}{Theorem}[section]

\newtheorem{lemma}[thm]{Lemma}

\newtheorem{preremark}[thm]{Remark}
 
\def\N{{\mathbb N}}

\def\R{{\mathbb R}} 
\def\a{{\alpha}} 
 
\newcommand{\CVD}{\hfill $\rule{2.6mm}{2.6mm}$} 
\newcommand{\PF}{\noindent{\bf Proof. }} 
 
\title[Quasilinear phase transitions]{ 
Rigidity results for some\\
boundary quasilinear phase transitions} 
 
\author[Y. Sire]{Yannick Sire} 
\author[E. Valdinoci]{Enrico Valdinoci}

\begin{document} 
\begin{abstract} 
We consider a quasilinear equation
given in the half-space, i.e. a so called
boundary reaction problem. Our concerns are a geometric Poincar\'e inequality 
and, as a byproduct of this inequality, a result on the symmetry of 
low-dimensional
bounded stable solutions, under some suitable assumptions on the nonlinearities. 
More precisely, we analyze the following boundary problem
$$
\left\{ 
\begin{matrix} 
-{\rm div}\, (a(x,|\nabla u|)\nabla u)+g(x,u)=0 \qquad 
{\mbox{ on $\R^n\times(0,+\infty)$}} 
\\ 
-a(x,|\nabla u|)u_x = f(u) 
\qquad{\mbox{ on $\R^n\times\{0\}$}}\end{matrix} 
\right.$$
under some natural assumptions on the diffusion coefficient
$a(x,|\nabla u|)$ and the nonlinearities $f$ and $g$. 

Here, $u=u(y,x)$, with~$y\in\R^n$ and~$x\in(0,+\infty)$.
This type of PDE can be seen as a nonlocal problem on the boundary 
$\partial \R^{n+1}_+$. The assumptions on $a(x,|\nabla u|)$ allow to treat in a unified way the $p-$laplacian and the minimal surface operators. 
\end{abstract} 
 
\maketitle 
\tableofcontents 
 
\bigskip\bigskip 
 
\bigskip 
 
\bigskip\bigskip 
 
\section{Introduction} 

The purpose of this paper is to give some geometric results on the following problem:
\begin{equation}\label{eq1-provv}
\left\{ 
\begin{matrix} 
-{\rm div}\, (a(x,|\nabla u|)\nabla u)+g(x,u)=0 \qquad 
{\mbox{ on $\R^n\times(0,+\infty)$}} 
\\ 
- a(x,|\nabla u|) u_x = f(u) 
\qquad{\mbox{ on $\R^n\times\{0\}$.}}\end{matrix} 
\right. 
\end{equation}
Here,
$u=u(y,x)$, with~$y\in\R^n$ and~$x\in(0,+\infty)$.
Equation \eqref{eq1-provv} is a boundary problem. This type of system is a model for nonlocal 
operators. For instance, when 
$g=0$ and $a(x,|\nabla u|)=x^\alpha $ with $\alpha \in (-1,1)$, it 
has been proved by~\cite{cafS}
that the Dirichlet-to-Neumann operator    
$$\Gamma: u|_{\partial \R^{n+1}_+} \mapsto -x^\alpha u_x|_{\partial \R^{n+1}_+}$$
is the fractional laplacian $(-\Delta)^{\frac{1-\alpha}{2}}$. 
In \cite{SV}, a symmetry result for bounded stable solutions of
semilinear equations involving this operator was given.

Unfortunately, a theory describing the boundary operator for problem \eqref{eq1-provv} is not yet
available. However, in virtue of the results
by~\cite{cafS}, one could
interpret the operator on the boundary as a nonlocal quasilinear operator.

In this paper, 
we develop
a geometric
analysis of the level sets
of stable solutions
of~\eqref{eq1-provv} and
we prove
a symmetry result
inspired by a conjecture
of De Giorgi~\cite{DeG}. 

We 
want to give a geometric insight of the phase  
transitions for equation~\eqref{eq1-provv}. Our goal 
is to give a geometric  
proof of the one-dimensional symmetry result
for boundary reactions in dimension~$n=2$, 
inspired by 
De Giorgi conjecture and 
in the spirit of the proof of 
Bernstein Theorem given in~\cite{giusti} and
applied in the case of boundary 
reactions in \cite{SV}.

We focus on problem \eqref{eq1-provv} under the following structural assumptions (denoted $(S)$):
\begin{itemize}

\item The function $a$ maps $(0,+\infty) \times
(0,+\infty)$ into $(0,+\infty)$ and
$$ \lim_{t\rightarrow 0^+} t a(.,t)=0.$$ 
\item The map $t \mapsto a(.,t)$ is $C^1(0,+\infty)$ and
\begin{equation}\label{BCOMEA}
t| a_t(x,t)| \leq C a(x,t)\end{equation}
for any~$x$, $t>0$,
for some constant $C>0$. 
\item The map $x \mapsto a(x,.)$ is in
$L^1((0,r))$, for any $r>0$ and bounded over all open sets compactly
contained in $\R^{n+1}_+$, i.e. for all $K \Subset 
\R^{n+1}_+$, 
there exists $\mu_1$,
$\mu_2>0$, possibly depending on $K$,
such that $\mu_1 \leq
a(x,t) \leq \mu_2$, for any $x\in K$ and for $0<t\leq M$.

Also, the function~$x \mapsto a(x,.)$ is an $A_2$-Muckenhoupt
weight, that is, there exists $\kappa>0$
such that
\begin{equation}\label{Muck}
\int_c^d a(x,t)\,dx \,
\int_c^d \frac{1}{a(x,t)}\,dx\,\le\,\kappa(d-c)^2 
\end{equation}
for any $d\ge c\ge 0$ and for all $0<t\leq M$.

\item The map $(0,+\infty)\ni x\mapsto g(x,0)$
belongs to $L^\infty((0,r))$ for any $r>0$.
Also,
for any $x>0$, the map $\R\ni u\mapsto g(x,u)$
is locally Lipschitz, and given any $R$, $M>0$
there exists $C>0$, possibly depending on $R$ and $M$
in such a way that
\begin{equation}\label{8ikeoqoqoqoo78}
\sup_{{0<x<R}\atop{|u|<M}}|g_u(x,u)|\le C.
\end{equation}
\item The function
$f$ is locally 
Lipschitz in
$\R$.  
\end{itemize}

Equation \eqref{eq1-provv} may be understood 
in the weak sense, namely supposing that $u\in 
L^\infty_{\rm loc}(\overline{\R^{n+1}_+})$, with  
\begin{equation}\label{hgasj7717177} 
a(x,|\nabla u|)|\nabla u|^2 \in L^1 (B_R^+) 
\end{equation} 
for any $R>0$, 
and that\footnote{Condition \eqref{hgasj7717177} 
is assumed here to make sense of \eqref{eq1}. 
We will see in the forthcoming Lemma \ref{Daf} that it is 
always uniformly fulfilled when $u$ is bounded and for a weight $a$ satisfying natural structural assumptions.

The structural
assumptions on $g$ may be easily
checked when $g(x,u)$ has the product-like
form of $g^{(1)}(x) g^{(2)}(u)$.} 
\begin{equation}\label{eq1} 
\int_{{\R^{n+1}_+}} 
a(x,|\nabla u|) \nabla u\cdot 
\nabla\xi+\int_{{\R^{n+1}_+}} g(x,u)\, \xi= 
\int_{\partial {\R^{n+1}_+}} 
f(u)\xi 
\end{equation} 
for any $\xi:B_R ^+\rightarrow \R$ which is bounded, locally
Lipschitz in the interior of
$\R^{n+1}_+$,
which 
vanishes on $\R^{n+1}_+\setminus B_R$ and such that
\begin{equation}\label{hgasj7717177-bis}
a(x,|\nabla u|)|\nabla\xi|^2
\in L^1 (B_R^+).\end{equation}

As usual,
we are using here the notation~$B_R^+:= B_R 
\cap\R^{n+1}_+$. 

Consider now the
map $\mathcal{B}: \R^+ \times \R^{n+1} \backslash \left \{ 0 \right \} 
\rightarrow {\rm{Mat}}
((n+1)\times (n+1))$ defined by
\begin{equation}
\label{BDE}
\mathcal{B}(x,\eta)_{ij}:=
a(x,|\eta|)\delta_{ij}+\frac{a_t(x,|\eta|)}{|\eta|} \eta_i \eta_j
\end{equation}
for any $1 \leq i,j \leq n+1$,
where $a_t$ stands for the derivative of $a(x,t)$ with respect to 
its second variable.

A direct computation gives 
\begin{equation}\label{6bis}
\frac{d}{d\varepsilon}
a(x,|\nabla u +\varepsilon
\nabla \varphi|)(\nabla u +\varepsilon
\nabla \varphi)\cdot \nabla \varphi |_{\varepsilon=0}=<\mathcal{B}
(x,\nabla u) \nabla \varphi,\nabla \varphi>
\end{equation}
for any smooth test function $\varphi$, any function $u$ with 
nonvanishing gradient and where $<,>$ stands for the canonical inner product in $\R^{n+1}$. 

Inspired by~\eqref{6bis}, it is tempting to say that~$u$ is 
stable
if
\begin{equation}\label{sta1} 
\begin{array}{c}
\int_{B_R^+} <\mathcal{B}(x,\nabla u)\nabla \xi, \nabla \xi>+
\int_{B_R^+} g_u(x,u)\xi^2 
-\int_{\partial B_R^+} 
f'(u)\xi^2\,\ge\,0 
\end{array}
\end{equation} 
for any $\xi$ as above. 
The above notion of
stability (sometimes 
also called semistability because of the large inequality) condition 
in~\eqref{sta1} appears naturally in the calculus 
of variations setting and it 
is usually related to minimization 
and monotonicity properties. 
In particular, \eqref{6bis} and~\eqref{sta1} 
state that the (formal) second variation 
of the energy functional associated 
to the equation has a sign (see, e.g.,~\cite{Moss, FCS, AAC} 
and Section~7 of~\cite{FSV} for further details). 

In our case, however, it is convenient to {\em relax} this definition
of stability. Namely, we say that~$u$ is stable
if~\eqref{sta1}
holds for any~$\xi$ of the form~$\xi:= |\nabla_y u|\phi$,
where~$\phi:\R^{n+1}\rightarrow \R$ is
Lipschitz and vanishes on $\R^{n+1}_+\setminus B_R$.

This relaxation of the stability definition is convenient for
our setting, since it makes possible to write~\eqref{sta1}
when $f$ is only locally Lipschitz and not necessarily differentiable.

Indeed, since the map~$y\mapsto u(y,x)$ will be taken to
be
locally Lipschitz (see~\eqref{LipA} below),
then so is the map~$y\mapsto f(u(y,x))$
and therefore
$$ f'(u)\xi^2 =\nabla_y \big(f(u)\big)\cdot \nabla_y u\,\phi^2$$
is well-defined almost everywhere, making sense of the last
term in~\eqref{sta1}.

The regularity assumption we take on~$u$ (see, in particular,~\eqref{hgasj7717177}
and~\eqref{SA3}) make the first term in~\eqref{sta1} well-posed too.
\bigskip 
 
The main results we prove are a geometric formula, 
of Poincar\'e-type, given in Theorem~\ref{POIN:TH}, 
and a symmetry result, given in Theorem~\ref{SYM:TH}. 
 
For our geometric result, we need to recall 
the following notation. Fixed $x>0$ and~$c\in\R$, we 
look at the level set 
$$ S:= \{  
y\in\R^n {\mbox{ s.t. }} 
u(y,x)=c 
\}.$$ 
We will consider the regular points of~$S$, 
that is, we define 
$$ L:=\{ y\in 
S 
{\mbox{ s.t. }} 
\nabla_y u(y,x)\neq 0 
\}.$$ 
Note that~$L$ depends on the~$x\in(0,+\infty)$ 
that we fixed at the beginning, though we do not keep 
explicit track of this in the notation. 
 
For any point $y\in L$, 
we let $\nabla_L$ to be the tangential gradient 
along~$L$, that is, for any~$y_o\in L$ 
and any~$G:\R^n\rightarrow\R$ smooth in the vicinity of~$y_o$, 
we set 
\begin{equation}\label{GR} \nabla_L G(y_o):= 
\nabla_y G(y_o)-\left(\nabla_y G(y_o)\cdot 
\frac{\nabla_y u(y_o,x)}{| 
\nabla_y u(y_o,x)|}\right) 
\frac{\nabla_y u(y_o,x)}{| 
\nabla_y u(y_o,x)|}.\end{equation}
Since~$L$ is a smooth manifold, in virtue of 
the Implicit Function Theorem (and of the standard 
elliptic 
regularity of $u$ 
apart from the boundary of $\R^{n+1}_+$), 
we can define 
the principal curvatures on it, denoted by 
$$\kappa_1(y,x),\dots, 
\kappa_{n-1}(y,x),$$ for any~$y\in L$. 
We will then define the total curvature 
$$ {\mathcal{K}}(y,x):=\sqrt{ 
\sum_{j=1}^{n-1} \big(\kappa_j (y,x)\big)^2 
}.$$ 
 
We also define 
$$ {\mathcal{R}}^{n+1}_+:=\{ 
(y,x)\in\R^n\times(0,+\infty){\mbox{ s.t. }} 
\nabla_y u(y,x)\neq 0 
\}.$$ 
 
With this notation, we can state 
our geometric formula: 

\begin{thm}\label{POIN:TH}
Assume that 
$u$ is a bounded and stable weak 
solution of~\eqref{eq1-provv} under assumptions $(S)$.

Assume furthermore that 
\begin{itemize}
\item For all $r>0$, 
\begin{equation}\label{LipA} |\nabla_y u|\in
{L^\infty(\overline{B_r^+})}. \end{equation}
\item For every $(y,x) \in B_R^+ \bigcap \left \{ \nabla u \neq 0 \right \}$, we have 
\begin{equation}\label{H1}
a(x,|\nabla u|) +\frac{a_t (x,|\nabla u|)}{|\nabla u|} u_x^2 \geq 0
\end{equation} 
and 
\begin{equation}\label{H2}
a(x,|\nabla u|) +\frac{a_t (x,|\nabla u|)}{|\nabla u|}
|\nabla_y u|^2 \geq \lambda(y,x) \ge0
\end{equation} 
for some $\lambda(y,x)$.
\end{itemize} 

Assume also the following regularity assumptions:
\begin{equation} \label{SA0}
\begin{split} 
&{\mbox{for almost any $x>0$, the map $\R^n 
\ni y\mapsto 
\nabla u(y,x)$}}\\ 
&{\mbox{ 
is in $W^{1,1}_{\rm loc}(\R^n, \R^n)$ ,
}}\end{split}\end{equation} 
\begin{equation}\label{SA3-provv} 
\begin{split} 
&{\mbox{the map $\R^{n+1}_+ 
\ni (y,x)\mapsto a(x,|\nabla u|) \sum_{j=1}^n 
\big( 
|\nabla u_{y_j}|^2+|u_{y_j}|^2 
\big)$}}\\ 
&{\mbox{is in $L^1(B_r^+)$, for any 
$r>0$
}}\end{split}\end{equation} 
and
\begin{equation}\label{SA3}
\begin{split} 
&{\mbox{the map $\R^{n+1}_+ 
\ni (y,x)\mapsto 
a(x,|\nabla u|) \big(  
|\nabla|\nabla_y u||^2+|\nabla_y u|^2 
\big) 
$}}\\ 
&{\mbox{is in
$L^1(B_r^+)$,
for any 
$r>0$. 
}}\end{split}
\end{equation} 
Then, 
for any $R>0$ and any~$\phi:\R^{n+1}\rightarrow \R$ which is 
Lipschitz and vanishes on $\R^{n+1}_+\setminus B_R$, we have that 
\begin{equation}\label{rtyua77a7a}
\begin{split}
&\int_{ {\mathcal{R}}^{n+1}_+} 
\,\phi^2 \left( 
a(x,|\nabla u|) {\mathcal{K}}^2 |\nabla_y u|^2+ 
\lambda(y,x) \big| 
\nabla_L |\nabla_y u| 
\big| ^2 
\right)\,\leq \\
&\, \int_{\mathcal{R}^{n+1}_+} |\nabla_y u|^2 <\mathcal{B}(x,\nabla u)\nabla \phi,\nabla \phi> .
\end{split}
\end{equation}
\end{thm} 

Assumption \eqref{LipA} is natural and it holds in particular in
the case $g := 0$, $a(x,t)=x^\a$ where $\a \in
(-1,1)$, as discussed in \cite{SV} (in many cases
of interest,
interior elliptic regularity then
ensures for free that $u$ is $C^1$
inside $\R^{n+1}_+$ and
$C^2$
as soon as the gradient does not vanish).
It is important to notice that
assumptions~\eqref{H1}
and~\eqref{H2} hold in the important
case of the $p-$laplace operator (i.e. $a(x,t)=t^{p-2}$, for~$p > 1$),
and in the case of the mean curvature operator (i.e. $a(x,t)=
\frac{1}{\sqrt{1+t^2}}$).

The regularity assumptions 
in~\eqref{SA0}, \eqref{SA3-provv}
and~\eqref{SA3}
are satisfied
in many cases of interest
(see, for instance Lemma~\ref{REGO} below).

The result in 
Theorem~\ref{POIN:TH} has been deeply inspired 
by the work of~\cite{SZarma, SZcrelle}, 
where related geometric inequalities have been first introduced
for the Allen-Cahn equation.
Further progress has been done
in~\cite{FAR, FSV} for reactions in the interior and
in~\cite{SV} for reactions on the boundary. 
 
The advantage of formula~\eqref{rtyua77a7a} is that 
one  
bounds 
tangential gradients and curvatures of level sets 
of stable solutions in terms of the gradient
of the solution
itself. 
That is, suitable 
geometric quantities of interest 
are controlled by an appropriate 
energy term. 
 
On the other hand, since the geometric formula bounds a 
weighted~$L^2$-norm of any 
test function~$\phi$ by a 
weighted~$L^2$-norm of its 
gradient, we may  
consider Theorem~\ref{POIN:TH} 
as a weighted {P}oincar\'e 
inequality. Again, the advantage of such a formula 
is that the weights have a neat geometric interpretation. 
See also~\cite{Ferrari} for further investigation
of Poincar\'e-type formulas.

The second result we present is a symmetry result 
in low dimension.  
 
\begin{thm}\label{SYM:TH}
Assume that $n=2$
and that the assumptions
in Theorem~\ref{POIN:TH} 
hold. Suppose also that~$\lambda(y,x)$
in~\eqref{H2}
is strictly positive almost everywhere.
Suppose also
that one of the following conditions~\eqref{g=0}
or~\eqref{g=+} hold, namely assume that
either for any $M>0$ 
\begin{equation}\label{g=0}
{\mbox{the map $(0,+\infty)\ni x\mapsto
\displaystyle\sup_{|u|\le M} |g (x,u)|$ is in $L^1((0,
+\infty))$}} 
\end{equation}
or that
\begin{equation}\label{g=+}
\inf_{{x\in\R^n}\atop{u\in\R}}g(x,u)\,u\,
\ge\,0.
\end{equation}

Assume that the diffusion coefficient $a(.,.)$ has a product structure
given by 
$$a(x,t)=\mu(x)\mathcal{A}(t),$$
where 
\begin{itemize}
\item the function $\mu$ is positive and such that
\begin{equation}\label{1.21bis}
\mu(x) \sim x^\alpha\end{equation}
for $\alpha \in (-1,1)$. 
\item One of the following two conditions is met:
either \begin{equation}\label{LB}
\mathcal{A} \in L^\infty (\R^+,\R^+)
\end{equation}
or 
\begin{equation}\label{PB}
\mathcal{A}(t) \sim t^{p-2}\end{equation}
with~$p \geq 1+\alpha$.
\end{itemize}
Then, there exist~$\omega:(0,+\infty)
\rightarrow {\rm S}^1$
and~$u_o: \R\times [0,+\infty)\rightarrow\R$
such that
$$ u(y,x)=u_o(\omega(x)\cdot y,x)$$
for any~$(y,x)\in {\R^{3}_+}$.
\end{thm} 

The paper~\cite{CSM} gave the first contribution
to symmetry result for boundary reaction PDEs.
In particular,~\cite{CSM} gave a result
analogous to Theorem~\ref{SYM:TH} when~$\mu:=1$, $g:=0$
and~$f\in C^{1,\beta}$.
In~\cite{SV}, a result analogous to Theorem~\ref{SYM:TH}
was given when~$a(x,t)=a(x)$, that is when~$a$ is
independent on the gradient term. In this sense, Theorem~\ref{SYM:TH}
extends the results of~\cite{CSM, SV}
to {\em quasilinear, possibly degenerate or singular, equations}
(in fact, when~$a(x,t):=x^\a$ and~$g:=0$, then~$\omega$ in Theorem~\ref{SYM:TH}
is constant, see~\cite{SV}).

We now discuss the assumptions of Theorem~\ref{SYM:TH}. First, 
the assumptions on~$\mathcal{A}$ are realized for
mean curvature operators, for which~$\mathcal{A}(t)=\frac{1}{\sqrt{
1+t^2}}$, which satisfies~\eqref{LB}
and for $p$-laplace operators, for which~$\mathcal{A}(t)=
t^{p-2}$,
when~$p\ge 1+\alpha$, which fulfills~\eqref{PB}.

The structural assumption on $\mu(x)$ 
is natural in the light of
the representation formula obtained in \cite{cafS}
which relates boundary reactions to fractional
operator (see also \cite{SV}): in this sense,
the operator studied here may be seen as a
quasilinear analogue of the fractional laplacian.

Theorem \ref{SYM:TH} 
asserts that, for any $x > 0$, the function $\R^2 
\ni y\mapsto 
u(y,x)$ depends only on one variable. Thus, 
Theorem \ref{SYM:TH} 
may be seen 
as the analogue of De Giorgi conjecture of~\cite{DeG} 
in dimension $n=2$ 
for equation \eqref{eq1-provv}. 

Condition \eqref{g=0}
is fulfilled by $g:=0$, or, more generally,
by $g:=g^{(1)}(x) g^{(2)}(u)$,
with $g^{(1)}$ summable over~$\R^+$
and~$g^{(2)}$ locally Lipschitz.
Also, condition~\eqref{g=+} is fulfilled by~$g:=u^{2\ell+1}$,
with~$\ell\in\N$.

When $u$ is not bounded,
the claim of Theorem \ref{SYM:TH} 
does not, in general, hold 
(a counterexample being~$a:=1$, $f:=0$,
$g:=0$ and~$
u(y_1,y_2,x):=y_1^2-y_2^2$).

Theorem~\ref{aux:P} below
will also provide a result, slightly
more general than Theorem \ref{SYM:TH},
which will be
valid for $n\ge 2$ and without
conditions~\eqref{g=0} or~\eqref{g=+},
under an additional energy assumption.
\bigskip

The rest of the paper is devoted to the proofs 
of Theorems~\ref{POIN:TH} and ~\ref{SYM:TH}. 
For this, some
preliminary energy estimate
will also be needed. 
 
\section{Some energy bounds}

This section is devoted to some preliminary energy estimate,
which are needed for the proof of
Theorem \ref{SYM:TH}.

Thus, throughout this section, the structural
assumptions
of Theorem \ref{SYM:TH} are in force.

We recall that
\begin{equation}
\label{9bis} a(x,|\nabla u|)\,u_x^2
\in L^1(B_R^+)
\end{equation}
for any $R>0$,
due to \eqref{hgasj7717177}.

We start with an elementary observation:

\begin{lemma}
There exists $C>0$ in such a way that
\begin{equation}\label{sii1kkkkj1j1}
\int_{B_{2R}^+ \setminus B_R^+} \mu(x)
\le CR^{n+1+\alpha}
\end{equation}
for any $R\ge1$ and $\alpha \in (-1,1)$.
\end{lemma}

\PF We have that
\begin{eqnarray*}
\int_{B_{2R}^+ \setminus B_R^+} \mu(x)
&\le&\int_{0}^{2R} \int_{B_{ 2R} } \mu(x)\,dy\,dx
\\ &\le& C_1 R^{n}
\int_{0}^{2R} \mu(x)\,dx\\
&\le& C_2 R^{n+1+\alpha},
\end{eqnarray*}
for suitable $C_1$, $C_2>0$,
due to~\eqref{1.21bis}.~\CVD
\medskip

Though not explicitly needed here, we would 
like to point out that the natural integrability 
condition in \eqref{hgasj7717177} 
holds uniformly for bounded solutions.
A byproduct of this gives an
energy estimate, which we will use in the proof of
Theorem \ref{SYM:TH}.
 
\begin{lemma}\label{Daf} 
For any $R>0$ there exists $C$, possibly 
depending on $R$, in such a way that
\begin{equation}\label{2.2bis} 
\| \mu(x) \mathcal{A}(|\nabla u|) |\nabla u|^2  \|_{L^1 (B_R^+)}\le C.
\end{equation}
Moreover, if 
\begin{itemize}
\item $n=2$, and
\item either \eqref{g=0} or~\eqref{g=+} holds,
\end{itemize} then there exists $C_o>0$
such that
\begin{equation}\label{AL}
\int_{B_R^+}\Big(
a(x,|\nabla u|) +|a_t(x,|\nabla u|)| \, |\nabla u|
\Big)\,|\nabla u|^2\,\le\, C_o\, R^2
\end{equation}
for any $R\ge 1$.
\end{lemma} 
 
\PF We focus on the proof of~\eqref{AL}, since~\eqref{2.2bis}
is a simple byproduct of the arguments we are going to perform.

The proof of Lemma~\ref{Daf}
consists in testing the weak formulation
in~\eqref{eq1}
with 
$\xi:=u
\tau ^\ell$ where $\tau$ is a cutoff
function
such that $0\le\tau\in C^\infty_0 (B_{2R})$, with $\tau=1$ 
in $B_{R}$ and $|\nabla \tau|\le 8/R$, with $R\ge1$. 
The parameter $\ell >1$ will be suitably chosen below.

Note that such a $\xi$ is admissible,
since~\eqref{hgasj7717177-bis}
follows from~\eqref{hgasj7717177}.

One then gets from \eqref{eq1} that 
\begin{eqnarray}\label{BAB}&&\nonumber 
\int_{{\R^{n+1}_+}}
a(x,|\nabla u|)\,
\big( |\nabla u |^2 \tau ^\ell +\ell  \tau^{\ell-1} u \nabla u \cdot \nabla \tau 
\big)+ 
\int_{{\R^{n+1}_+}} g(x,u) u \tau^\ell\\&&\qquad=
\int_{\R^n} f(u) u\tau^\ell. 
\end{eqnarray} 

We now distinguish the case in which~\eqref{LB} holds
from the case in which~\eqref{PB} holds.

If~\eqref{LB} holds,
we take $\ell =2$. 
Thus, by Cauchy-Schwarz
inequality, we deduce from~\eqref{BAB} that
\begin{eqnarray*}
\int_{\R_+^{n+1}}\mu(x)\mathcal{A}(|\nabla u|)\,|\nabla u|^2\tau^2
\le \frac 12 \int_{\R_+^{n+1}}\mu(x)\mathcal{A}(|\nabla u|)\,|\nabla u|^2\tau^2
\\
+
C_* \Big(
\int_{\R_+^{n+1}}
\mu(x)\mathcal{A}(|\nabla u|)|\nabla \tau|^2+\int_{\R^{n}}|f(u)|\,|u|\,\tau^2
\Big)-
\int_{\R_+^{n+1}}
g(x,u)\,u\,\tau^2,
\end{eqnarray*}
for a suitable constant $C_*>0$.

This, recalling \eqref{BCOMEA},
\eqref{g=0}, \eqref{g=+}, \eqref{LB} and~\eqref{sii1kkkkj1j1},
plainly gives~\eqref{AL}.

If, on the other hand,~\eqref{PB} holds,
we take $\ell =p$. Therefore,we have 
\begin{equation*}
\int_{\R_+^{n+1}}\mu(x)\mathcal{A}(|\nabla u|)\,|\nabla u|^2\tau^p
\sim \int_{\R_+^{n+1}}\mu(x)\,|\nabla u|^p\tau^p.
\end{equation*}
Recalling~\eqref{BAB} and
using \eqref{g=0}, \eqref{g=+}, \eqref{sii1kkkkj1j1}, one has 
\begin{equation*}
\int_{\R_+^{n+1}}\mu(x)\,|\nabla u|^p\tau^p \leq C
\Big\{
\int_{\R_+^{n+1}}\mu(x)\,|\nabla u|^{p-1} |\nabla \tau | \tau^{p-1} +R^n
\Big\}.  
\end{equation*}
Thus, by Young inequality, we conclude that
\begin{eqnarray*}
\int_{\R_+^{n+1}}\mu(x) |\nabla u|^p\tau^p
\le C\Big \{ \varepsilon \int_{\R_+^{n+1}}
\Big \{ \mu(x)^{1/q}|\nabla u|^{p-1} \tau^{p-1} \Big \} ^q +\\ 
C_\varepsilon \int_{\R_+^{n+1}}\mu(x)|\nabla \tau
|^p +R^n \Big \}
\end{eqnarray*}
for some $\varepsilon>0$ and $q=\frac{p}{p-1}$. 

Making use of~\eqref{sii1kkkkj1j1},
this leads to 

\begin{eqnarray*}
\int_{\R_+^{n+1}}\mu(x) |\nabla u|^p\tau^p
\le C\Big \{\int_{B_{2R}}\frac{\mu(x)}{R^p}+R^n \Big \}
\leq  C(R^{n+1+\alpha-p} +R^n). 
\end{eqnarray*}

This gives the desired result as soon
as $p \geq 1+\alpha.$~\CVD\medskip  

\section{The Poincar\'e-type formula: proof of Theorem \ref{POIN:TH}} 
 
This section is devoted to the proof
of the geometric formula in Theorem \ref{POIN:TH}.
As we will see throughout the proof,
the assumptions in Theorem \ref{POIN:TH} are natural and quite general. 

Besides few technicalities, the 
proof of Theorem \ref{POIN:TH} consists
in plugging the right test function in stability 
condition~\eqref{sta1} and in using the linearization 
of \eqref{eq1-provv} to get rid of the unpleasant terms. 
Following are the rigorous details of the proof. 
   
By~\eqref{BDE}, we have that 
\begin{equation}\label{26bis}
\begin{split}
&\int_{\mathcal{R}^{n+1}_+ } a(x,|\nabla u|)\nabla u \cdot \Psi_{y_j}=
\\&
-\int_{\mathcal{R}^{n+1}_+ } a(x,|\nabla u|) \nabla u_{y_j} \cdot \nabla \Psi + a_t
(x,|\nabla u|) \frac{\nabla u \cdot \nabla u_{y_j}}{|\nabla u|} \nabla u \cdot \Psi=\\
&
=-\int_{\mathcal{R}^{n+1}_+ } <\mathcal{B}(x,\nabla u) \nabla u_{y_j}, \Psi >. 
\end{split}\end{equation}
for any~$j=1,\dots, n$ and any~$\Psi\in C^\infty 
(\R^{n+1}_+, \R^n)$ supported in~$B_R$. 
 
Making use of~\eqref{eq1} and~\eqref{26bis} with~$\Psi:=\nabla\psi$,
we conclude that 
\begin{equation}\label{a711aa} 
\begin{split} 
&\int_{\R^{n+1}_+} g_u(x,u) u_{y_j} \psi - 
\int_{\R^n} f'(u) u_{y_j} \psi=
\\
&\int_{\R^{n+1}_+} (g(x,u))_{y_j} \psi - \int_{\R^n} (f(u))_{y_j}\psi=\\
&-\int_{\R^{n+1}_+} g(x,u)\psi_{y_j} + \int_{\R^n} f(u)\psi_{y_j}=
\\&-\int_{\mathcal{R}^{n+1}_+  } <\mathcal{B}(x,\nabla u) \nabla u_{y_j}, \nabla \psi >
\end{split} 
\end{equation} 
for any~$j=1,\dots, n$ and any~$\psi\in 
C^\infty 
(\R^{n+1}_+)$ supported in~$B_R$. 
 
A density argument (see, e.g.,
Lemma~3.4, Theorem~2.4 and~(2.9) in~\cite{CPSC})
via~\eqref{BCOMEA}
and~\eqref{SA3-provv}, implies that~\eqref{a711aa} 
holds for~$\psi:=u_{y_j} \phi^2$, 
where~$\phi$ is 
as in the statement of 
Theorem~\ref{POIN:TH}, therefore 
\begin{equation}\label{3.2bis}
\begin{split}
& 0= \int_{ {B}_R^+}
<\mathcal{B}(x,\nabla u) \nabla u_{y_j},
\nabla u_{y_j} > \phi^2+ <\mathcal{B}(x,\nabla u) \nabla u_{y_j}, \nabla \phi^2> u_{y_j}+\\
& \int_{B_R^+} g_u(x,u) u_{y_j}^2 \phi^2 -\int_{\partial B_R^+} f'(u) u_{y_j}^2 \phi^2.
\end{split}
\end{equation}

Let now~$r$, $\rho>0$ and~$P\in \R^{n+1}_+$ be such
that~$B_{r+\rho}(P)\subset \R^{n+1}_+$.
We consider~$\gamma$ to be either~$|\nabla_y u|$
or~$u_{y_j}$. In force of~\eqref{SA3-provv}
and~\eqref{SA3}, we see that~$\gamma$
is in~$W^{1,2}( B_r(P))$, and so in~$W^{1,1}_{\rm loc}
(B_r (P))$.

Thus, by Stampacchia Theorem (see, e.g., Theorem~6.19
in~\cite{LOSS}), $\nabla \gamma=0$ for almost
any~$(y,x)\in B_r(P)$ such that~$\gamma(y)=0$.

Hence, since~$P$, $r$ and~$\rho$ can be
chosen arbitrarily,
we have that
\begin{equation}\label{ASTA}{\mbox{
$\nabla |\nabla_y u| =0=\nabla u_{y_j}$
for almost every~$(y,x)$ such that~$\nabla_y u(y,x)=0$.}}\end{equation}

By~\eqref{3.2bis} and~\eqref{ASTA}, we obtain
\begin{equation*}
\begin{split} 
& 0= \int_{\mathcal{B}_R^+} 
<\mathcal{B}(x,\nabla u) \nabla u_{y_j},
\nabla u_{y_j} > \phi^2+ <\mathcal{B}(x,\nabla u) \nabla u_{y_j}, \nabla \phi^2> u_{y_j}+\\
& \int_{B_R^+} g_u(x,u) u_{y_j}^2 \phi^2 -\int_{\partial B_R^+} f'(u) u_{y_j}^2 \phi^2. 
\end{split} 
\end{equation*} 
where $\mathcal{B}_R^+=B_R^+ \bigcap \mathcal{R}^{n+1}_+.$
We now sum over $j=1,...,n$ to get (dropping, for short,
the dependences of $\mathcal{B}$) and
we obtain
\begin{equation} \label{78iddudududuudaa}
\begin{split} 
&-\int_{\mathcal{B}_R^+} \sum_{j=1}^n <\mathcal{B}\nabla u_{y_j}, \nabla u_{y_j} >
\phi^2- \frac{1}{2}<\mathcal{B} \nabla |\nabla_y u|^2, \nabla \phi^2>=\\
& \int_{B_R^+} g_u(x,u) |\nabla_y u|^2 \phi^2 -\int_{\partial B_R^+} f'(u) |\nabla_y u|^2 \phi^2. 
\end{split} 
\end{equation} 

Now, we make use of~\eqref{sta1} 
by taking~$\xi:=|\nabla_y u|\phi$ 
(this choice was also performed 
in~\cite{SZarma, SZcrelle, FAR, FSV,SV}; 
note that \eqref{LipA} 
and~\eqref{SA3} imply \eqref{hgasj7717177-bis}
and so they
make it possible to 
use here such a test function). We thus obtain 
\begin{equation*}
\begin{split}
& 0 \leq \int_{\mathcal{B}_R^+} <\mathcal{B}
\nabla |\nabla_y u|, \nabla |\nabla_y u|> \phi^2 
+ <\mathcal{B} \nabla \phi, \nabla \phi> |\nabla_y u|^2+\\
& 2 <\mathcal{B}\nabla |\nabla_y u| ,
\nabla \phi> |\nabla_y u| \phi 
+ g_u(x,u) |\nabla_y u|^2 \phi ^2 -\int_{\partial B_R^+} f'(u) |\nabla_y u|\phi^2,
\end{split} 
\end{equation*}
where~\eqref{ASTA} has been used once more.

This and~\eqref{78iddudududuudaa}
imply that 
\begin{equation}\begin{split} 
\label{s88818181} 
&0 \leq \int_{\mathcal{B}_R^+} <\mathcal{B} \nabla |\nabla_y u|, \nabla |\nabla_y u|> \phi^2 +<\mathcal{B} \nabla \phi, \nabla \phi> |\nabla_y u|^2 \\
&-\sum_{j=1}^n <\mathcal{B}\nabla u_{y_j}, \nabla u_{y_j} > \phi^2. 
\end{split}\end{equation} 

By using \eqref{BDE} and~\eqref{s88818181}, we are lead
to the following inequality
\begin{equation}\label{31jkl}
\begin{split}
& 0 \leq \int_\mathcal{B_R^+} a(x,|\nabla u|) \phi^2 \Big [|\nabla |\nabla_y u||^2
-\sum_{j=1}^n  |\nabla u_{y_j}|^2\Big ]+\\&\quad
< \mathcal{B} \nabla \phi, \nabla \phi> |\nabla_y 
u|^2 +\\
&\quad\quad
\frac{a_t(x,|\nabla u|) \phi^2}{|\nabla u|} \Big [(\nabla u \cdot \nabla |\nabla_y 
u|)^2 -\sum_{j=1}^n (\nabla u \cdot \nabla u_{y_j})^2 \Big ]
.\end{split}
\end{equation}

We denote
$$
\mathcal{H}_*:=
-(\partial_x |\nabla_y u|)^2+\sum_{j=1}^n u_{xy_j}^2 ,$$
$$\mathcal{H}_1:=|\nabla |\nabla_y u||^2-\sum_{j=1}^n |\nabla u_{y_j}|^2$$
$${\mbox{and }}\quad
\mathcal{H}_2=:(\nabla u \cdot \nabla |\nabla_y u |)^2-\sum_{j=1}^n (\nabla u 
\cdot \nabla u_{y_j})^2.$$
 
We have that
\begin{equation}
\label{A1}\begin{split}
& \mathcal{H}_2= (u_x \partial_x |\nabla_y u|)^2 -\sum_{j=1}^n (u_x u_{xy_j})^2+ 
(\nabla_y u \cdot \nabla_y |\nabla_y u|)^2-\sum_{j=1}^n (\nabla_y u \cdot \nabla_y u_{y_j})^2
\\ &\qquad=
-u_x^2 \mathcal{H}_*
+(\nabla_y u \cdot \nabla_y |\nabla_y u|)^2-\sum_{j=1}^n (\nabla_y u \cdot \nabla_y u_{y_j})^2
,\end{split}\end{equation}
where we have just separated the $x$ and $y$ variables. 

Also, {f}rom~\eqref{GR},
\begin{equation}\label{3.5bis}
|\nabla_L G|^2=|\nabla_y G|^2-
\left(\nabla_y G \cdot \frac{\nabla_y u}{|\nabla_y u|}
\right)^2,\end{equation}
for any smooth function~$G:\R^n\rightarrow\R$.

Hence, making use of~\eqref{3.5bis}
with~$G:=|\nabla_y u|$, we obtain that,
on~${{\mathcal{R}}^{n+1}_+}$,
\begin{equation}
\label{A2}
\begin{split}
&(\nabla_y u \cdot \nabla_y |\nabla_y u|)^2-\sum_{j=1}^n (\nabla_y u \cdot \nabla_y u_{y_j})^2=\\
&|\nabla_y u|^2 \Big [ \Big(
\frac{\nabla_y u}{|\nabla_y u|} \cdot \nabla_y |\nabla_y 
u|\Big)^2-\sum_{j=1}^n \Big(
\frac{\nabla_y u}{|\nabla_y u|} \cdot \nabla_y u_{y_j}
\Big)^2 \Big ]=\\
&|\nabla_y u|^2 \Big [ |\nabla_y |\nabla_y u| |^2 -|\nabla_L |\nabla_y u ||^2 -\sum_{j=1}^n
\Big(\frac{\nabla_y u}{|\nabla_y u|} \cdot \nabla_y u_{y_j}
\Big)^2\Big ]=\\
&-|\nabla_y u|^2 |\nabla_L |\nabla_y u ||^2.
\end{split}
\end{equation}

By a differential geometry formula
obtained in~\cite{SZarma, SZcrelle}
(see also
equation~(2.10) in~\cite{FSV}), we have, on~${{\mathcal{R}}^{n+1}_+}$, 
\begin{equation}\label{A3}
\mathcal{H}_1=-\mathcal{H}_* -(\mathcal{K}^2 |\nabla_y u|^2+|\nabla_L |\nabla_y u||^2).
\end{equation}

As a consequence of~\eqref{A1},
\eqref{A2} and~\eqref{A3}, we obtain that~\eqref{31jkl}
may be written in the following form: 
\begin{eqnarray} 
\label{yuiooaoo} 
\nonumber 
0 \leq \int_{{{\mathcal{R}}^{n+1}_+}} a(x,|\nabla u|) \phi^2 \Big (-\mathcal{H}_* -(\mathcal{K}^2|\nabla_y u|^2 +|\nabla_L |\nabla_y u||^2)\Big )\\
+\frac{a_t(x,|\nabla u|) \phi^2}{|\nabla u|} \Big (-u_x^2 \mathcal{H}_* -|\nabla_y u|^2 |\nabla_L |\nabla_y u||^2 \Big )+\\
<\mathcal{B}\nabla \phi,\nabla \phi > |\nabla_y u |^2.\nonumber 
\end{eqnarray} 
We now note that, on~${{\mathcal{R}}^{n+1}_+}$, by Cauchy-Schwarz inequality, we have 
$\mathcal{H}_* \geq 0.$ 

This,~\eqref{yuiooaoo} and assumptions \eqref{H1}-\eqref{H2} complete the proof of 
Theorem~\ref{POIN:TH}.~\CVD 
 
\section{The symmetry result: proof of Theorem \ref{SYM:TH}} 
 
As in~\cite{FSV, SV},
the strategy for proving Theorem \ref{SYM:TH} 
is to test the geometric formula of 
Theorem~\ref{POIN:TH} against an appropriate capacity-type 
function to make the left  
hand side vanish. 
This would give that the curvature of the level sets for fixed~$x>0$ 
vanishes and so that these level sets are flat, as desired 
(for this, 
the vanishing of the tangential gradient term is also 
useful to take care of the possible 
plateaus of~$u$, where 
the level sets are not smooth manifold).

As described in the assumptions of Theorem \ref{SYM:TH}, we will
take some structure for the weight $a(x,|\nabla u|)$ (in fact,
such assumptions might be further weakened, paying the price
of additional 
technicalities in the proofs).

Some preparation is needed for the proof 
of Theorem \ref{SYM:TH}. 
Indeed, Theorem \ref{SYM:TH} will 
follow from the subsequent Theorem~\ref{aux:P}, which 
is valid for any dimension~$n$ and without the restriction
in either \eqref{g=0}
or \eqref{g=+}. 
 
We will use the notation~$X:=(y,x)$ for points in~$\R^{n+1}_+$.  
 
Given~$\rho_1\le\rho_2$, we also define 
$${\mathcal{A}}_{\rho_1,\rho_2}:=\{ 
X\in\R^{n+1}_+{\mbox{ s.t. }}|X|\in [\rho_1,\rho_2] 
\}.$$ 
 
\begin{lemma}\label{tatay} 
Let~$R>0$ 
and~$h:B_R^+\rightarrow\R$ be a nonnegative 
measurable function.  
 
For any~$\rho\in (0,R)$,
let 
$$ \eta(\rho):=\int_{B^+_{\rho}} h.$$ 
Then, 
$$\int_{{\mathcal{A}}_{\sqrt R, R}}\frac{h(X)}{|X|^2}\,dX 
\leq 2\int_{\sqrt R}^R t^{-3}\eta(t)\,dt+\frac{\eta(R)}{R^2}. 
$$ 
\end{lemma} 
 
For the proof
of Lemma~\ref{tatay}, see Lemma~10 in~\cite{SV}.
 
\begin{thm}\label{aux:P} 
Let $u$ be as requested in Theorem \ref{POIN:TH}.
Assume furthermore that  
there exists~$C_o\geq 1$ in such a way that 
\begin{equation}\label{en:bound} 
\int_{B^+_R}  \Big(
a(x,|\nabla u|)+|a_t(x,|\nabla u|)|\, |\nabla u|
\Big)|\nabla u|^2\le C_o\,  
R^2\end{equation} 
for any~$R\ge C_o$. 
 
Then there exist~$\omega:(0,+\infty)
\rightarrow {\rm S}^1$
and~$u_o: \R\times [0,+\infty)\rightarrow\R$
such that
$$ u(y,x)=u_o(\omega(x)\cdot y,x)$$ 
for any~$(y,x)\in\R^{n+1}_+$. 
\end{thm} 
 
\PF {F}rom Lemma~\ref{tatay} applied here with
$$h(X):=
\Big(a(x,|\nabla u(X)|)+|a_t(x,|\nabla u(X)|)|\, |\nabla u(X)|
\Big) |\nabla 
u(X)|^2$$ and \eqref{en:bound}, we obtain
\begin{equation}\label{7s77s88} 
\begin{split}&
\int_{{\mathcal{A}}_{\sqrt R, R}}\frac{\Big(a(x,|\nabla u(X)|)+
|a_t(x,|\nabla u(X)|)|\, |\nabla 
u(X)|\Big) |\nabla u(X)|^2 
}{ |X|^2}\\ &\qquad\leq C_1\log R 
\end{split}\end{equation} 
for a suitable~$C_1$, 
as long as~$R$ is large enough. 
 
Now we define 
$$ \phi_R(X):=\left\{ 
\begin{matrix} 
\log R & {\mbox{ if $|X|\le \sqrt R$,}}\\ 
2\log\big( R/|X|\big)\Big)  
& {\mbox{ if $\sqrt R<|X|< R$,}} 
\\ 
0 & {\mbox{ if $|X|\ge R$}} 
\end{matrix} 
\right.$$ 
and we observe that 
\begin{equation}\label{4.2bis}
|\nabla\phi_R|\leq \frac{C_2\,\chi_{ 
{\mathcal{A}}_{\sqrt R, R} 
}}{|X|},\end{equation}
for a suitable~$C_2>0$. 

From~\eqref{BDE}
and Cauchy-Schwarz
inequality, we have that,
for any~$w\in \R^{n+1}$,
\begin{equation}\label{4.2tris}
|<\mathcal{B}(x,\nabla u)
w, w >| \leq \Big \{a(x,|\nabla u|)+|a_t(x,|\nabla u|)||\nabla u| \Big \} |
w|^2.\end{equation}
Thus, plugging~$\phi_R$ inside the geometric 
inequality of Theorem~\ref{POIN:TH}, we obtain 
\begin{eqnarray*} 
&& (\log R)^2\int_{B^+_{\sqrt{R}}\bigcap 
{\mathcal{R}}^{n+1}_+ 
} 
\left( a(x,|\nabla u|)
{\mathcal{K}}^2 |\nabla_y u|^2+ \lambda(y,x)
\big| 
\nabla_L |\nabla_y u| 
\big| ^2 
\right)\\&&\qquad\qquad\,\leq\,C_3 
\int_{ 
{\mathcal{A}}_{\sqrt R, R} 
}\frac{ 
\Big(
a(x,|\nabla u|)+|a_t(x,|\nabla u|)||\nabla u|
\Big)\,
|\nabla_y u|^2}{|X|^2} 
\end{eqnarray*} 
for large~$R$, thanks to~\eqref{4.2bis} and~\eqref{4.2tris}.
 
By dividing by~$(\log R)^2$, 
employing~\eqref{7s77s88} 
and taking~$R$ arbitrarily large, we conclude
that~${\mathcal{K}}$ and~$\big| 
\nabla_L |\nabla_y u| 
\big|$ vanish identically
on~${\mathcal{R}}^{n+1}_+$. 
 
Then, the desired result follows 
by Lemma~2.11 of~\cite{FSV} (applied to the function~$y\mapsto 
u(y,x)$, for any fixed~$x>0$).~\CVD 
\medskip 

We now complete the proof of 
Theorem~\ref{SYM:TH}.
We observe that, under the assumptions of
Theorem~\ref{SYM:TH},
estimate \eqref{en:bound}
holds, thanks to \eqref{AL}.
Consequently, the hypotheses of Theorem~\ref{SYM:TH} 
imply the ones of Theorem~\ref{aux:P}, 
from which the claim in Theorem~\ref{SYM:TH} follows.~\CVD 

\section{Further comments on assumptions
\eqref{SA0},
\eqref{SA3-provv} and~\eqref{SA3}}

Having completed the proof of the main results,
in this section we would like to remark
that assumptions~\eqref{SA0},
\eqref{SA3-provv} and~\eqref{SA3}
are quite natural in many cases of interest.

For instance, we assume in this section that
the structural hypotheses on~$a(x,t)$ in Theorem~\ref{SYM:TH}
and the bound in~\eqref{LipA}
hold true.

For simplicity,
we also suppose that~$u$ is~$C^2_{\rm loc} (\R^{n+1}_+)$
(this is the case, for instance, 
of mean curvature type operators or of~$p-$laplace
operators if~$\nabla u$ does not vanish). The purpose
of this section is then to show that
conditions~\eqref{SA0},
\eqref{SA3-provv} and~\eqref{SA3}
are satisfied in this case.

\begin{lemma}\label{Nuovo} 
We have
$$\mu(x) \mathcal{A}(|\nabla u|)|\nabla u_{y_j}|^2 \in L^1(B_R^+)$$ 
for every $R>0$.  
\end{lemma} 
 
\PF 
Given $|\eta|<1$, $\eta\ne 0$, we consider the incremental
quotient
$$ u_\eta(y,x):= \frac{u(y_1,\dots,y_j+\eta,\dots,y_n,x)-
u(y_1,\dots,y_j,\dots,y_n,x)}{\eta}.$$
Since $f$ is locally Lipschitz,
\begin{equation}\label{Al1}
[f(u)]_\eta \le C,
\end{equation}
for some $C>0$, due to \eqref{LipA}.

Analogously, from \eqref{8ikeoqoqoqoo78}
and \eqref{LipA},
for any $R>0$ there exists $C_R>0$ such that
\begin{equation}\label{Al2}
[g(x,u)]_\eta \le C_R
\end{equation}
for any $x\in(0,R)$.

Let now $\xi$ be as requested in~\eqref{eq1}.
Then,~\eqref{eq1} gives that
\begin{eqnarray}\label{5.2bis}
\nonumber&& 
\int_{{\R^{n+1}_+}}\big[
\mu(x) \mathcal{A}(|\nabla u|)\nabla u_{\eta}\cdot 
\nabla\xi+\big(g(x,u)\big)_\eta \,\xi\big]- 
\int_{\partial {\R^{n+1}_+}} \big[
f(u)\big]_\eta\xi\\ 
&=& 
-\int_{{\R^{n+1}_+}}\big[
\mu(x) \mathcal{A}(|\nabla u|)\nabla u\cdot 
\nabla\xi_{-\eta}+g(x,u) \,\xi_{-\eta} 
\big]+
\int_{\partial {\R^{n+1}_+}} 
f(u) \xi_{-\eta}\\ 
&=&0. \nonumber
\end{eqnarray} 

We concentrate on the case when~\eqref{LB} holds
(the case in which~\eqref{PB} holds is then an
easy
modification, analogous to the one performed
in the proof of Lemma~\ref{Daf}).

We consider a smooth cutoff function
$\tau$ such that $0\le\tau\in C^\infty_0 (B_{R+1})$, with $\tau=1$ 
in $B_{R}$ and $|\nabla \tau|\le 2$. 
Taking $\xi:=u_{\eta}\tau^2$ in~\eqref{5.2bis}, one gets 
\begin{equation}\label{0a8h1hmclakk} 
\begin{split}
& 2 \int_{{\R^{n+1}_+}}
\mu(x) \mathcal{A}(|\nabla u|)\tau u_{\eta} \nabla u_{\eta
}\cdot \nabla \tau \\&\quad 
+\int_{{\R^{n+1}_+}} 
\mu(x) \mathcal{A}(|\nabla u|)\tau^2 |\nabla u_{\eta}|^2
+\int_{{\R^{n+1}_+}}  \big(g(x,u)\big)_\eta
u_{\eta}\,\tau^2\\&\quad\quad
=\int_{\partial {\R^{n+1}_+}} 
\big(f(u)\big)_\eta\,u_{\eta} \tau^2.
\end{split}\end{equation}
We remark that the above choice of~$\xi$
is admissible, since
\eqref{hgasj7717177-bis}
follows from~\eqref{LipA} and \eqref{9bis}.
  
Now, by  Cauchy-Schwarz
inequality, we have
\begin{equation*}
\begin{split} 
&\int_{{\R^{n+1}_+}}
\mu(x) \mathcal{A}(\nabla u|)\tau u_{\eta} \nabla u_{\eta}\cdot \nabla \tau \geq
-\frac{\varepsilon}{2} \int_{{\R^{n+1}_+}}
\mu(x) \mathcal{A}(|\nabla u|)\tau^2 |\nabla u_{\eta}|^2\\
& -\frac{1}{2\varepsilon} \int_{{\R^{n+1}_+}}
\mu(x) \mathcal{A}(|\nabla u|)u_{\eta}^2|\nabla \tau|^2
\end{split}
\end{equation*}
for any $\varepsilon>0$. 

Therefore, by choosing $\varepsilon$ suitably small,~\eqref{0a8h1hmclakk}
reads
\begin{eqnarray*}
&&
\int_{{\R^{n+1}_+}}
\mu(x) \mathcal{A}(|\nabla u|)\tau^2 |\nabla u_{\eta}|^2
\\
&\le& C\,
\Big[
\int_{B_{R+1}^+}
\mu(x) \mathcal{A}(|\nabla u|)u_{\eta}^2+
\int_{B_{R+1}^+}  \big| \big(g(x,u)\big)_\eta
u_{\eta}\big|
\\
&&\quad+
\int_{\{|y|\le R\}\times\{x=0\}} \big|
\big(f(u)\big)_\eta
u_{\eta}
\big|\Big].
\end{eqnarray*}
for some $C>0$.

{F}rom~\eqref{LipA}, 
\eqref{LB},
\eqref{Al1} and \eqref{Al2},
we thus control
$$\int_{B_R^+}
\mu(x) \mathcal{A}(|\nabla u|)\tau^2 |\nabla u_{\eta}|^2
$$
uniformly in $\eta$.

By sending $\eta\rightarrow 0$
and using Fatou Lemma,
we obtain the desired claim.~\CVD 
\medskip

Following is the regularity needed for
some subsequent computations.

\begin{lemma} \label{REGO}
Conditions~\eqref{SA0},
\eqref{SA3-provv} and~\eqref{SA3}
are satisfied.\end{lemma}

The proof is omitted, since it is analogous
to the one of Lemma~7 in~\cite{SV}.

\section*{Acknowledgments} 
 
YS would like to thank
the hospitality of Universit\`a di
Roma Tor Vergata, where part of this work has been 
done.  
 
EV has been partially supported by~{\em MIUR 
Me\-to\-di va\-ria\-zio\-na\-li ed equa\-zio\-ni 
dif\-fe\-ren\-zia\-li non\-li\-nea\-ri}.

\bibliographystyle{alpha} 
\bibliography{bibliofile} 
 
\vfill
 
{\em YS} --  
Universit\'e Aix-Marseille 3, Paul C\'ezanne -- 
LATP -- 
Marseille, France. 
 
{\tt sire@cmi.univ-mrs.fr} 
\medskip 
 
{\em EV} -- 
Universit\`a di Roma Tor Vergata -- 
Dipartimento di Matematica -- 
I-00133 Rome, Italy. 
 
{\tt valdinoci@mat.uniroma2.it}
\end{document}